\documentclass[10pt]{amsart}
\usepackage{latexsym,amsmath,amssymb,amscd}
\def\today{\ifcase \month \or
   January \or February \or March \or April \or
   May \or June \or July \or August \or
   September \or October \or November \or December \fi
   \space\number\day , \number\year}

\begin{document}

\makeatletter
\@addtoreset{figure}{section}
\def\thefigure{\thesection.\@arabic\c@figure}
\def\fps@figure{h,t}
\@addtoreset{table}{bsection}

\def\thetable{\thesection.\@arabic\c@table}
\def\fps@table{h, t}
\@addtoreset{equation}{section}
\def\theequation{
\arabic{equation}}
\makeatother

\newcommand{\bfi}{\bfseries\itshape}

\def\theoremname{Theorem}
\def\propositionname{Proposition}
\def\corollaryname{Corollary}
\def\lemmaname{Lemma}
\def\remarkname{Remark}
\def\stepname{Step}
\def\definitionname{Definition}
\def\exercisename{Exercise}
\def\examplename{Example}
\def\examplesname{Examples}
\def\problemname{Problem}
\def\problemsname{Problems}
\def\proofname{Proof}

\def\@thmcounter#1{\noexpand\arabic{#1}}
\def\@thmcountersep{}
\def\@begintheorem#1#2{\it \trivlist \item[\hskip 
\labelsep{\bf #1\ #2.\quad}]}
\def\@opargbegintheorem#1#2#3{\it \trivlist
      \item[\hskip \labelsep{\bf #1\ #2.\quad{\rm #3}}]}
\makeatother
\newtheorem{theo}{\theoremname}[section]
\newtheorem{propo}[theo]{\propositionname}
\newtheorem{coro}[theo]{\corollaryname}
\newtheorem{lemm}[theo]{\lemmaname}

\newenvironment{theorem}{\begin{theo}\it}{\end{theo}}
\newenvironment{proposition}{\begin{propo}\it}{\end{propo}}
\newenvironment{corollary}{\begin{coro}\it}{\end{coro}}
\newenvironment{lemma}{\begin{lemm}\it}{\end{lemm}}

\newtheorem{rem}[theo]{\remarkname}
\newenvironment{remark}{\begin{rem}\rm}{\end{rem}}

\newtheorem{defin}[theo]{\definitionname} 
\newenvironment{definition}{\begin{defin}\rm}{\end{defin}}

\newtheorem{notat}[theo]{Notation} 
\newenvironment{notation}{\begin{notat}\rm}{\end{notat}}

\newtheorem{ex}[theo]{\examplename}
\newenvironment{example}{\begin{ex}\rm}{\end{ex}}

\newtheorem{exs}[theo]{\examplesname}
\newenvironment{examples}{\begin{exs}\rm}{\end{exs}}

\newtheorem{conj}[theo]{\conjecturename}
\newenvironment{conjecture}{\begin{conj}\rm}{\end{conj}}

\newtheorem{pr}[theo]{\problemname}
\newenvironment{problem}{\begin{pr}\rm}{\end{pr}}

\newcommand{\todo}[1]{\vspace{5 mm}\par \noindent
\framebox{\begin{minipage}[c]{0.85 \textwidth}
\tt #1 \end{minipage}}\vspace{5 mm}\par}

\newcommand{\1}{{\bf 1}}

\newcommand{\hotimes}{\widehat\otimes}

\newcommand{\Ci}{{\mathcal C}^\infty}
\newcommand{\comp}{\circ}
\newcommand{\D}{\text{\bf D}}
\newcommand{\de}{{\rm d}}
\newcommand{\ev}{{\rm ev}}
\newcommand{\id}{{\rm id}}
\newcommand{\ie}{{\rm i}}
\newcommand{\End}{{\rm End}\,}
\newcommand{\Hom}{{\rm Hom}\,}
\newcommand{\Ker}{{\rm Ker}\,}
\newcommand{\lf}{{\rm l}}
\newcommand{\Ran}{{\rm Ran}\,}
\newcommand{\rad}{{\rm rad}\,}
\newcommand{\spann}{{\rm span}}
\newcommand{\Tr}{{\rm Tr}\,}
\newcommand{\trile}{\trianglelefteq}

\newcommand{\ad}{\mathop{{\rm ad}}\nolimits}
\newcommand{\Exp}{\mathop{{\rm Exp}}\nolimits}
\newcommand{\Log}{\mathop{{\rm Log}}\nolimits}
\newcommand{\GL}{\mathop{{\rm GL}}\nolimits}
\newcommand{\PGL}{\mathop{{\rm PGL}}\nolimits}

\newcommand{\G}{{\rm G}}
\newcommand{\U}{{\rm U}}
\newcommand{\VB}{{\rm VB}}

\newcommand{\Ac}{{\mathcal A}}
\newcommand{\Bc}{{\mathcal B}}
\newcommand{\Cc}{{\mathcal C}}
\newcommand{\Dc}{{\mathcal D}}
\newcommand{\Gc}{{\mathcal G}}
\newcommand{\Hc}{{\mathcal H}}
\newcommand{\Ic}{{\mathcal I}}
\newcommand{\Kc}{{\mathcal K}}
\newcommand{\Lc}{{\mathcal L}}
\newcommand{\Nc}{{\mathcal N}}
\newcommand{\Oc}{{\mathcal O}}
\newcommand{\Pc}{{\mathcal P}}
\newcommand{\Qc}{{\mathcal Q}}
\newcommand{\Sc}{{\mathcal S}}
\newcommand{\Xc}{{\mathcal X}}
\newcommand{\Zc}{{\mathcal Z}}

\newcommand{\Gg}{{\mathfrak G}}
\newcommand{\Hg}{{\mathfrak H}}
\newcommand{\Jg}{{\mathfrak J}}
\newcommand{\Sg}{{\mathfrak S}}

\newcommand{\g}{{\mathfrak g}}
\newcommand{\h}{{\mathfrak h}}
\newcommand{\fj}{{\mathfrak j}}
\newcommand{\s}{{\mathfrak s}}
\newcommand{\fr}{{\mathfrak r}}

\newcommand{\N}{{\mathbb N}}
\newcommand{\Z}{{\mathbb Z}}
\newcommand{\R}{{\mathbb R}\,}
\newcommand{\C}{{\mathbb C}\,}
\newcommand{\K}{{\mathbb K}}
\newcommand{\F}{{\mathbb F}}
\newcommand{\Q}{{\mathbb Q}}
\renewcommand{\H}{{\mathbb H}}

\renewcommand{\:}{\colon}
\newcommand{\0}{{\bf 0}}
\newcommand{\striang}{{\frak{sut}}}
\newcommand{\triang}{{\frak ut}}

\newcommand{\subeq}{\subseteq}
\newcommand{\supeq}{\supseteq}
\newcommand{\ctr}{{}_\rightharpoonup}
\newcommand{\into}{\hookrightarrow}
\newcommand{\eps}{\varepsilon}

\def\onto{\to\mskip-14mu\to} 

\renewcommand{\hat}{\widehat}
\renewcommand{\tilde}{\widetilde}

\newcommand{\Rarrow}{\Rightarrow}
\newcommand{\nin}{\noindent} 
\newcommand{\oline}{\overline}
\newcommand{\Fa}{{\rm Fa}}
\newcommand{\Larrow}{\Leftarrow}
\newcommand{\la}{\langle}
\newcommand{\ra}{\rangle}
\newcommand{\Mon}{{\rm Mon}}
\newcommand{\up}{\mathop{\uparrow}}
\newcommand{\down}{\mathop{\downarrow}}
\newcommand{\res}{\vert}
\newcommand{\sdir}{\times_{sdir}}

\renewcommand{\L}{\mathop{\bf L{}}\nolimits}

\pagestyle{myheadings}
\markboth{\sl CIA}{\sl CIA}


\makeatletter
\title{Finite-dimensional Lie subalgebras of algebras with continuous inversion} 
\author{Daniel~Belti\c t\u a} 
\author{Karl-Hermann~Neeb}
\address{Institute of Mathematics ``Simion
Stoilow'' of the Romanian Academy, 
P.O. Box 1-764, RO-014700 Bucharest, Romania}
\email{Daniel.Beltita@imar.ro}
\address{Department of Mathematics, Darmstadt University of Technology, 
Schlossgartenstrasse 7, D-64289 Darmstadt, Germany}
\email{neeb@mathematik.tu-darmstadt.de}
\date{January 21st, 2008}

\begin{abstract}
We investigate the finite-dimensional Lie groups 
whose points are separated by the continuous homomorphisms 
into groups of invertible elements of locally convex algebras 
with continuous inversion that satisfy an appropriate completeness condition. 
We find that these are precisely the linear Lie groups, 
that is, 
the Lie groups which can be faithfully represented as matrix groups. 
Our method relies on proving that certain finite-dimensional Lie subalgebras 
of algebras with continuous inversion commute modulo the Jacobson radical. 


\end{abstract}
\makeatother
\maketitle

\section{Introduction}

In the present paper we investigate the finite-dimensional Lie groups 
whose points are separated by the continuous homomorphisms 
into groups of invertible elements of locally convex algebras 
with continuous inversion. 
We find that these are precisely the linear Lie groups, 
that is, 
the Lie groups which can be faithfully represented as matrix groups 
(see Theorem~\ref{final} below). 
And the structure of the linear groups is fairly well understood; 
see for instance Hochschild's book~\cite{Ho65}. 

This problem is motivated by recent developments in the Lie theory of infinite-dimensional 
Lie groups modeled on locally convex spaces (cf.\ \cite{Mil83}, \cite{Bel06}, \cite{GN06}). 
In this context, the unit groups of continuous 
inverse algebras are the prototypical ``linear Lie groups'' 
(\cite{Gl02}), and it is a natural 
question whether the notion of ``linearity'' in this general context 
determines a larger class of finite-dimensional Lie groups than 
the Lie groups of matrices. 

The unital Banach algebras provide examples of locally convex algebras 
with continuous inversion, 
and in this special case we recover the characterization obtained 
in the paper~\cite{LV94} for the Lie groups 
whose uniformly continuous representations separate the points. 
However, it is not clear to us if the approach used in \cite{LV94}  
can be extended to non-normable algebras. 
For one thing, one of the key tools used in the aforementioned paper 
was a version of Lie's Theorem concerning 
existence of weights for infinite-dimensional representations  
of solvable Lie algebras. 
And such a version is available as yet only for representations 
by Banach space operators 
(see \cite{GL73} and \cite{BS01}). 
It is worth mentioning at this point that there exist large classes of
algebras with continuous inversion which are not Banach algebras,
for instance algebras of germs of holomorphic functions or
algebras of smooth matrix-valued functions on compact manifolds;
see Examples VIII.3 in \cite{Ne06} for more details and 
additional examples.

Thus an alternative approach is needed when working with general 
locally convex algebras with continuous inversion. 
For this purpose we prove that the main result of \cite{Ti87} 
actually holds true far beyond the setting of Banach algebras, 
where it was originally discovered. 
Loosely speaking, we show that if $\g$ is a finite-dimensional 
complex solvable Lie subalgebra of an algebra with continuous inversion 
that satisfies an appropriate completeness condition, 
then the closed unital associative subalgebra generated by $\g$ 
is commutative modulo its (Jacobson) radical; 
see Theorem~\ref{turovskii} below, 
which is the main technical result we need here. 
We obtain it by a method inspired from Turovskii's proof under the version 
exposed in~\cite{BS01}. 
However we emphasize that the present paper can be read 
independently of that book, inasmuch as 
several of the key tools developed in~\cite{BS01} 
---notably the Kleinecke-Sirokov Theorem and 
the spectral theory for several non-commuting variables---  
are not readily available beyond the setting of normable algebras, 
so that we now replaced them by arguments of a different type. 
And we strove to take advantage of this situation in order to make the present 
paper fairly self-contained. 

In Section~2 we set up the necessary preliminaries 
on algebras with continuous inversion 
and on nilpotent elements in Lie subalgebras of associative algebras. 
Section~3 concerns spectra of commutators, and includes in particular 
a version of the Kleinecke-Sirokov Theorem holding for not necessarily normable algebras,  
as well as a version of Rosenblum's theorem from \cite{Ro56} suitable for our purposes. 
In Section~4 we prove the theorem on commutativity modulo the radical 
(Theorem~\ref{turovskii}). 
Finally, in Section~5 we obtain the main result of the present paper, 
namely the characterization of the finite-dimensional Lie groups 
whose points are separated by the homomorphisms into 
groups of invertible elements of FC-complete algebras with continuous inversion.

\section{Preliminaries}

\subsection*{Preliminaries on algebras with continuous inversion}

\medbreak

\begin{notation}\label{bounded}
For an arbitrary unital complex associative algebra $\Ac$ we shall 
use the following notation: 
\begin{itemize}
\item[$\bullet$] $\Ac^\times=\{a\in\Ac\mid(\exists a^{-1}\in\Ac)\quad 
aa^{-1}=a^{-1}a=\1\}$;
\item[$\bullet$] the {\it spectrum} of any $a\in\Ac$ is 
$\sigma_{\Ac}(a) = \sigma(a)
=\{\lambda\in{\mathbb C}\mid \lambda\1-a\not\in\Ac^\times\}$.
\item[$\bullet$] the {\it spectral radius} of any $a\in\Ac$ is 
$r_{\Ac}(a)=\sup\{\vert\lambda\vert\mid\lambda\in\sigma(a)\}\in[0,\infty]$;
\item[$\bullet$] the {\it center} 
$\Zc_{\Ac}=\{a\in\Ac\mid(\forall b\in\Ac)\quad ab=ba\}$; 
\item[$\bullet$] $\Nc_{\Ac}=\{a\in\Ac\mid(\exists N\in{\mathbb N})\quad 
a^N=0\}$; 
\item[$\bullet$] $\Qc_{\Ac}=\{a\in\Ac\mid\sigma(a)=\{0\}\} = \{ a \in \Ac \mid \1+ \C a 
\subseteq \Ac^\times\}$; 
\item[$\bullet$] the {\it radical} $\rad\Ac=
\{a\in\Ac\mid(\forall b\in\Ac)\quad \1-ab\in\Ac^\times\} \subeq \Qc_\Ac$
\end{itemize}
Moreover, for any complex vector space $\Xc$ we denote by $\Lc(\Xc)$ 
the set of all linear maps from $\Xc$ into itself.
\end{notation}

\begin{definition} 
A {\it continuous inverse algebra} (CIA for short) is 
a Hausdorff locally convex unital algebra $\Ac$ whose unit group $\Ac^\times$ is 
open and for which the inversion map $\Ac^\times \to \Ac, a \mapsto a^{-1}$ 
is continuous. 

For any element $a$ of a complex CIA $\Ac$ the spectrum $\sigma(a)$ is a non-empty 
compact subset of $\C$. Here the boundedness of the spectrum follows from 
$$ r_{\Ac}(a) = (\sup \{ r > 0 \mid |z| \leq r \Rarrow \1 + z a \in \Ac^\times\})^{-1} $$
and its closedness from the openness of the unit group $\Ac^\times$. 

If, in addition, $\Ac$ is complete, then the same arguments 
as for Banach algebras lead to a holomorphic functional calculus 
(\cite{Wae67}, \cite{Gl02}). 
Since completeness is in general not inherited 
by quotients (\cite{Koe69}, \S 31.6), it is natural to consider for 
CIAs the weaker condition that they are {\it FC-complete} in the sense 
that they are closed under holomorphic functional calculus. This means 
that for $a \in A$, any open neighborhood $U$ of $\sigma(a)$, 
each holomorphic function $f \in \Oc(U)$ and any contour $\Gamma$ 
around $\sigma(a)$ in $U$, the integral 
$$ f(a) := \frac{1}{2\pi i} \oint_\Gamma f(\zeta) (a - \zeta \1)^{-1}\, d\zeta, $$
which defines an element of the completion of $\Ac$, actually exists in $\Ac$. 
\end{definition}

\begin{remark}\label{radical}
\normalfont 
It is known that in every unital complex associative algebra $\Ac$ 
the radical $\rad\Ac$ is equal to the intersection of 
all the maximal left ideals of $\Ac$, 
and it is equal to the intersection of 
all the maximal right ideals of $\Ac$ as well. 
(See for instance Theorem~3.53 in Chapter~1 of \cite{He93}.) 
In particular $\rad\Ac$ is a two-sided ideal of $\Ac$, 
and we have $\rad(\Ac/\rad\Ac)=\{0\}$. 
(See e.g., Theorem~3.63 in Chapter~1 of \cite{He93}.)

If $\Ac$ is a CIA, then the fact that its unit group is open implies that 
every maximal left ideal of $\Ac$ is closed, and this is the case 
with the right ideals as well (see e.g., Subsection~2.3 in Chapter~II of \cite{Wae67}, 
or Corollary~3.9 in Chapter~2 of~\cite{Co68}), 
hence $\rad\Ac$ is a closed two-sided ideal of~$\Ac$. 
Then the quotient algebra $\Ac/\rad\Ac$ is in turn 
a locally convex unital algebra with continuous inversion. 
(See e.g., Subsection~2.2 in Chapter~II of~\cite{Wae67} 
on quotients by closed ideals, 
or Prop.~3.14 in Chapter~2 of~\cite{Co68}.)
\end{remark}

\begin{lemma}\label{radii}
If $\Ac$ is a commutative complex CIA, 
then the spectral radius $r_{\Ac}$ is a continuous submultiplicative seminorm 
and, in particular, $\rad \Ac = \Qc_{\Ac}$
\end{lemma}

\begin{proof} (cf.\ \cite{Bi04}, Th.~1.7) The spectrum 
$\hat \Ac := \Hom(\Ac, \C)$ of $\Ac$ is a compact Hausdorff space and 
the Gelfand transform 
$\Gc \: \Ac \to C(\hat A), a \mapsto \hat a, \hat a(\chi) := \chi(a)$
is a continuous homomorphism of complex unital locally convex algebras, satisfying 
$\hat a(\hat \Ac) = \sigma(a)$ and hence 
$\|\hat a\|_\infty = r_{\Ac}(a)$. From that it follows that 
for each quasi-nilpotent element $a \in \Ac$ and $b \in \Ac$ we have 
$r_{\Ac}(ab) \leq r_{\Ac}(a)r_{\Ac}(b)=0$, so that $\1 - ab \in \Ac^\times$, 
which leads to $a \in \rad \Ac$. 
This implies that $\rad \Ac \subeq \Qc_\Ac \subeq \rad \Ac$, 
which leads to the asserted equality. 
\end{proof}

An early reference for our Lemma~\ref{radii} is  for instance \cite{Wae67}.
Another proof follows by 
Lemma~II.9, Proposition~II.3, and Corollary~III.9 in \cite{KOO98}.  

\begin{lemma}\label{subalgebras}
Let $\Ac$ be a complex CIA. 
\begin{itemize}
\item[(1)] Each closed unital subalgebra $\Sc$ of $\Ac$ is a CIA. 
If $\Ac$ is FC-complete, then so is $\Sc$. 
\item[(2)] If $\Sc \leq \Ac$ is a maximal commutative subalgebra, 
then $\Sc$ is unital, closed and equispectral, 
i.e., for each $x \in \Sc$ we have 
$\sigma_{\Ac}(x) = \sigma_{\Sc}(x)$. 
\item[(3)] For each closed ideal $\Ic \trile \Ac$ the quotient algebra 
$\Qc := \Ac/\Ic$ is a CIA. If, in addition, $\Ac$ is FC-complete, and the quotient map 
$q \: \Ac \to \Qc$ satisfies $\sigma_\Ac(a) = \sigma_{\Qc}(q(a))$ for each 
$a \in A$, then $\Qc$ is also FC-complete. 
\end{itemize}
\end{lemma}

\begin{proof} (1) (See also Remarque~4.6 and Subsection~4.3 in Chapter~2 of \cite{Wae67}.)  

If $a \in \Sc$ satisfies $r_{\Ac}(a) < 1$, then $\1 - a$ is invertible and 
the Neumann series $\sum_{n = 0}^\infty a^n$ converges to $(\1 - a)^{-1}$ 
(\cite{Gl02}). Since 
$\Sc$ is closed, $(\1 -a)^{-1} \in \Sc$, so that $\Sc^\times$ is a neighborhood of 
$\1$ in $\Sc$, hence open. 
The continuity of the inversion in $\Sc$ follows from the corresponding 
property of $\Ac$.

If, in addition, $\Ac$ is FC-complete and $s \in \Sc$, then 
$\sigma_\Sc(s) \supeq \sigma_\Ac(s)$, so that each contour around 
$\sigma_{\Sc}(s)$ also surrounds $\sigma_\Ac(s)$. 
Now for any holomorphic function $f$ on a neighborhood 
of $\sigma_{\Sc}(s)$ the integral 
$$ f(s) :=  \frac{1}{2\pi i}\oint_\Gamma f(\zeta) (s - \zeta \1)^{-1}\, d\zeta $$
defines an element of $\Ac$. Since $(s - \zeta \1)^{-1} \in \Sc$ for each $\zeta$, 
the closedness of $\Sc$ yields $f(s) \in \Sc$. This shows that $\Sc$ is FC-complete. 

(2) Since closures of commutative subalgebras are commutative subalgebras, 
the maximality of $\Sc$ implies its closedness. It trivially implies $\1 \in \Sc$. 
Moreover, for each $x \in \Sc$ we have 
$(x - \lambda \1)^{-1} \in \Sc$ whenever $x - \lambda \1$ is invertible 
in $\Ac$ because $(x - \lambda \1)^{-1}$ commutes with all elements of $\Sc$. 
This leads to $r_{\Ac}(x) = r_{\Sc}(x)$. 

(3) Since $q(\Ac^\times) \subeq \Qc^\times$ is an open subset of $\Qc$, 
the unit group $\Qc^\times$ is open. The continuity of the inversion of $\Qc$ 
follows from its continuity in $\1$ and 
the continuity of the map $\Ac^\times \to \Qc^\times, a \mapsto q(a)^{-1} = q(a^{-1})$, 
because $q$ is an open map. 

For each $a \in A$ we have $\sigma_{\Ac}(a)= \sigma_\Qc(q(a))$. 
For any open neighborhood of $\sigma_{\Qc}(q(a))$, 
each holomorphic function $f \in \Oc(U)$, and any contour around $\sigma(a)$ in $U$ 
the FC-completeness implies that the integral 
$$ f(a) :=  \frac{1}{2\pi i}\oint_\Gamma f(\zeta) (a - \zeta \1)^{-1}\, d\zeta $$
exists in $\Ac$. We conclude that the integral 
$$ q(f(a)) 
=  \frac{1}{2\pi i}\oint_\Gamma f(\zeta) q(a - \zeta \1)^{-1}\, d\zeta
=  \frac{1}{2\pi i}\oint_\Gamma f(\zeta) \big(q(a) - \zeta \1\big)^{-1}\, d\zeta = f(q(a)) $$
exists in $\Qc$. Hence $\Qc$ is FC-complete. 
\end{proof}

\subsection*{Algebraic preliminaries}
In this subsection we turn to the purely algebraic part of the preliminaries 
we need for our main results. In particular, we 
prove a suitable generalization of Theorem~2 in \S 28 of \cite{BS01}.

\begin{proposition}\label{theorem2}
Let $\Xc$ be a complex vector space, $\g$ a finite-dimensional 
solvable Lie subalgebra of $\Lc(\Xc)$, 
and $\h$ a Cartan subalgebra of $\g$. 
Denote by $\g^\alpha$, $\alpha\in R$, the family of root spaces of 
$\g$ corresponding to the set $R$ of non-zero roots, 
and assume that $\g^\alpha$ consists of nilpotent elements for every $\alpha\in R$. 
Then $\Nc_{\Lc(\Xc)}\cap\g$ is an ideal of $\g$. 
\end{proposition}

\begin{proof}
The proof has two steps. 

{\bf Step 1:} If $\g$ is a nilpotent Lie algebra 
then the desired conclusion follows by precisely the same reasoning 
as in Step 1 of the proof of Theorem~2 in \S 28 of \cite{BS01}. 
Specifically, we shall prove the desired conclusion by induction on 
$\dim\g$. 
The assertion is clear if $\dim\g=1$. 

Now assume that $\dim\g>1$ and $\Nc_{\Lc(\Xc)}\cap\g\ne\{0\}$. 
Since nilpotent elements are always polynomially central 
(Definition~2 in \S 16 of \cite{BS01}), 
it follows by Theorem~1 in \S 18 of \cite{BS01}, applied with $I = \g$ 
(see also \cite{Sa96}), that 
there exists $Y\in\Nc_{\Lc(\Xc)}\cap\g$ 
such that $Y\ne0$ and $[Y,T]=0$ for all $T\in\g$. 
Let $m\ge1$ with $Y^{m-1}\ne0=Y^m$, 
and denote $\Xc_k=\Ker(Y^k)$ for $k=0,\dots,m$. 
Since $Y$ belongs to the center of $\g$, it follows that 
$\{0\}=\Xc_0\subseteq\cdots\subseteq\Xc_m=\Xc$ is a nest of 
invariant subspaces for $\g$. 
In particular, there exist representations 
$$\rho_k\colon\g\to\Lc(\Xc_k/\Xc_{k-1}),\quad 
\rho_k(T)(x+\Xc_{k-1})=Tx+\Xc_{k-1}
\text{\ \ for\ all\ \ } T\in\g\text{\ and\ }x\in\Xc_k,$$
for $k=0,\dots,m$. 
It is easy to see that we have 
\begin{equation}\label{intersection}
\Nc_{\Lc(\Xc)}\cap\g
=\bigcap_{k=0}^n\rho_k^{-1}(\Nc_{\Lc(\Xc_k/\Xc_{k-1})})
=\bigcap_{k=0}^n\rho_k^{-1}(\Nc_{\Lc(\Xc_k/\Xc_{k-1})}\cap\rho_k(\g)),
\end{equation}
that is, $T\in\g$ is nilpotent if and only if each 
$\rho_k(T)\in\Lc(\Xc_k/\Xc_{k-1})$ is nilpotent for $k=0,\dots,m$. 

On the other hand $0\ne Y\in\bigcap\limits_{k=0}^m\Ker\rho_k$, 
hence for each $k$ we have $\dim\rho_k(\g)<\dim\g$. 
Then the induction hypothesis shows that 
$\Nc_{\Lc(\Xc_k/\Xc_{k-1})}\cap\rho_k(\g)$ is an ideal of $\rho_k(\g)$. 
Now \eqref{intersection} implies that $\Nc_{\Lc(\Xc)}\cap\g$ is an 
ideal of $\g$ since pull-backs and intersections of ideals 
are again ideals. 

{\bf Step 2:} We now proceed with the proof in the general case. 
Denote $\g_+:=\bigoplus\limits_{\alpha\in R}\g^\alpha$ 
so that the generalized root space decomposition of $\g$ 
with respect to the Cartan subalgebra $\h$ 
leads to the Fitting decomposition $\g=\h\dotplus\g_+$. 
Next denote 
$\g_0=\{T\in\g\mid\ad_{\g}T\colon\g\to\g\text{ is nilpotent}\}$. 
Since $\g$ is a solvable Lie algebra, it easily follows by the Lie theorem 
on simultaneous triangularization of solvable Lie algebras of matrices 
that $\g_0$ is an ideal of $\g$. 
Moreover, $\g_0$ is a nilpotent Lie algebra and we have
$$\bigcup\limits_{\alpha\in R}\g^\alpha
\subseteq\Nc_{\Lc(\Xc)}\cap\g\subseteq\g_0.$$
In particular Step 1 of the proof 
shows that $\Nc_{\Lc(\Xc)}\cap\g$ is an ideal of $\g_0$. 
To prove that $\Nc_{\Lc(\Xc)}\cap\g$ is an ideal of $\g$ 
it remains to check that 
$[\h,\Nc_{\Lc(\Xc)}\cap\g]\subseteq\Nc_{\Lc(\Xc)}\cap\g$. 

To this end note that since we have the inclusion of vector subspaces 
$\g_+\subseteq\Nc_{\Lc(\Xc)}\cap\g$ 
and $\g=\h\dotplus\g_+$ 
it follows that 
$\Nc_{\Lc(\Xc)}\cap\g=(\Nc_{\Lc(\Xc)}\cap\h)\dotplus\Cc_{\h}$. 
Since $\g=\h\dotplus\g_+$ it then follows that 
\begin{equation}\label{three}
[\h,\Nc_{\Lc(\Xc)}\cap\g]\subseteq
[\h,\Nc_{\Lc(\Xc)}\cap\h]+[\h,\g_+].\end{equation}
Now $[\h,\Nc_{\Lc(\Xc)}\cap\h]\subseteq\Nc_{\Lc(\Xc)}\cap\h$ 
by Step~1 of the proof, 
since $\h$ is a nilpotent Lie algebra. 
Moreover, $[\h,\g_+]\subseteq\g_+\subseteq\Nc_{\Lc(\Xc)}\cap\g$, 
where the latter inclusion follows by the hypothesis since 
we saw that $\Nc_{\Lc(\Xc)}\cap\g$ is a vector space. 
Thus \eqref{three} shows that 
$[\g,\Nc_{\Lc(\Xc)}\cap\g]\subseteq\Nc_{\Lc(\Xc)}\cap\g$, 
and the proof is complete. 
\end{proof}

\begin{lemma}\label{lemma5}
Assume that $\Ac$ is a complex unital associative algebra 
and $\g$ a finite-dimensional Lie subalgebra of $\Ac$ 
such that $\g\subseteq\Nc_{\Ac}$. 
Then there exists an integer $m\ge1$ such that 
$a_1\cdots a_m=0$ for all $a_1,\dots,a_m\in\g$. 
\end{lemma}

\begin{proof} Since $\ad a$ is a nilpotent operator on $\Ac$ for each $a \in \g$, 
Engel's Theorem implies that $\g$ is a nilpotent Lie algebra. 
In view of the Poincar\'e-Birkhoff-Witt Theorem,  
the unital associative subalgebra $\Ac(\g)$ generated by $\g$ 
is finite-dimensional. Then Lie's Theorem, applied to the 
left regular representation of $\g$ on $\Ac(\g)$ implies that 
the associative subalgebra of $\Ac$ generated by $\g$ consists of 
nilpotent elements. 
\end{proof}

\section{Spectra of commutators}

\begin{proposition}\label{fine}
Let $\Ac$ be a complex FC-complete CIA. 
Assume that $a,b,c\in\Ac$ satisfy 
$[a,b] = ab-ba=c$, $ac=ca$, and $bc=cb$. 
Then $\sigma_{\Ac}(c)=\{0\}$. 
\end{proposition}

\begin{proof} Since $\Ac$ is FC-complete, it has an exponential 
function $\exp \: \Ac \to \Ac^\times, x \mapsto e^x$, defined by 
the holomorphic functional calculus. 
Define 
$$f\colon{\mathbb C}\to\Ac,\quad f(t)=e^{ta}be^{-ta}.$$
Then 
$f$ is holomorphic and we have 
$f(0)=b$ and $f'(t)=ae^{ta}be^{-ta}-e^{ta}bae^{-ta}=e^{ta}ce^{-ta}=c$ 
for all $t\in{\mathbb C}$, because of the assumption $ac=ca$. 
This implies that 
$e^{ta}be^{-ta}=b+tc$ for each $t \in \C$, and hence that 
\begin{equation*}
\sigma_{\Ac}(b)=\sigma_{\Ac}(b+tc) \quad \mbox{ for each} \quad t \in \C.
\end{equation*}

As $b$ and $c$ commute, 
$$ |t|\sigma(c) = \sigma(tc) = \sigma(tc + b - b) 
\subeq \sigma(tc+b) - \sigma(b)
= \sigma(b) - \sigma(b). $$
Since $\sigma(b)$ is bounded, we obtain for $|t| \to \infty$ 
the inclusion $\sigma(c) \subeq \{0\}$. The lemma now follows 
from the non-emptyness of the spectrum. 
\end{proof}

\begin{lemma} \label{center} 
For each complex FC-complete CIA $\Ac$ we have $\Zc_\Ac\cap\Qc_\Ac\subseteq\rad\Ac$. 
\end{lemma}

\begin{proof} Let $c\in\Zc_\Ac\cap\Qc_\Ac$ and $b\in\Ac$ arbitrary. 
We have to show that $\1-bc\in\Ac^\times$. 
To this end, let $\Ac_0$ be a maximal commutative subalgebra of $\Ac$ 
containing both $b$ and $c$. Then $\Ac_0$ is a closed unital subalgebra 
of $\Ac$, hence an FC-complete CIA with 
$r_{\Ac}(a) = r_{\Ac_0}(a)$ (Lemma~\ref{subalgebras}). 
In particular $r_{\Ac}(bc)=r_{\Ac_0}(bc)\le r_{\Ac_0}(b)r_{\Ac_0}(c)=0$, 
where the inequality follows by Lemma~\ref{radii}. 
Thus $\sigma_{\Ac}(bc)=\{0\}$, and then $\1-bc\in\Ac^\times$, as desired.
\end{proof}

\begin{proposition}\label{rosenblum}
Let $\Ac$ be a complex FC-complete CIA. 
Then for all $a_1,a_2\in\Ac$ the operator 
$$\Delta\colon\Ac\to\Ac,\quad x\mapsto a_1x-xa_2,$$
satisfies 
$\sigma_{\Lc(\Ac)}(\Delta)\subseteq\sigma_{\Ac}(a_1)-\sigma_{\Ac}(a_2)$.
\end{proposition}

\begin{proof}
The method of proof used in Section~3 of \cite{Ro56} 
works in the present setting as well. 
Specifically, let $\lambda\in{\mathbb C}$ 
such that $\lambda\not\in\sigma_{\Ac}(a_1)-\sigma_{\Ac}(a_2)$. 
We are going to prove that 
$\lambda\not\in\sigma_{\Lc(\Ac)}(\Delta)$, either.
Since both $\sigma_{\Ac}(a_1)$ and $\sigma_{\Ac}(a_2)$ are compact subsets 
of $\C$, there exist two open subsets 
$U_1$ and $U_2$ of ${\mathbb C}$ such that 
$\sigma_{\Ac}(a_j)\subseteq U_j$ for $j=1,2$ 
and $\overline{U}_1\cap(\lambda+\overline{U}_2)=\emptyset$. 

Now let $\Gamma\subseteq U_2$ be a piecewise smooth contour 
surrounding $\sigma_{\Ac}(a_2)$. 
Then for every $z\in\Gamma$ we have $\lambda+z\not\in\sigma_{\Ac}(a_1)$, 
hence we can define a linear map from $\Ac$ into itself by 
$$T\colon\Ac\to\Ac,\quad 
Tx=\frac{1}{2\pi\ie}\oint\limits_{\Gamma}
((\lambda+z)\1-a_1)^{-1}x(z\1-a_2)^{-1}\de z.$$
We now have 
\begin{eqnarray*}
T(\lambda \1- \Delta)x 
&=& T(\lambda x - a_1 x + x a_2) \\
&=& \frac{1}{2\pi i} \oint_\Gamma ((\lambda + z\1) - a_1)^{-1} 
(\lambda x - a_1 x + x a_2) (z\1 - a_2)^{-1} \de z  \\
&=& \frac{1}{2\pi i} \oint_\Gamma ((\lambda + z\1) - a_1)^{-1} 
\big(((\lambda +z)\1 - a_1)x + x(a_2 - z \1)\big) (z\1 - a_2)^{-1} \de z  \\
&=& \frac{1}{2\pi i} \oint_\Gamma x (z\1 - a_2)^{-1} \de z  
- \frac{1}{2\pi i} \oint_\Gamma ((\lambda + z\1) - a_1)^{-1} x\de z  \\
&=& x\frac{1}{2\pi i} \oint_\Gamma (z\1 - a_2)^{-1} \de z  
- \frac{1}{2\pi i} \oint_\Gamma ((\lambda + z\1) - a_1)^{-1} \de z \cdot x  \\
&=& x\cdot \1 - \0 \cdot x = x 
\end{eqnarray*}
and 
\begin{eqnarray*}
(\lambda \1- \Delta)Tx 
&=& (\lambda Tx - a_1 Tx + Tx a_2) \\
&=& \frac{1}{2\pi i} \oint_\Gamma (\lambda\1 - a_1)((\lambda + z\1) - a_1)^{-1} x
 (z\1 - a_2)^{-1} \de z  \\
&& + \frac{1}{2\pi i} \oint_\Gamma ((\lambda + z\1) - a_1)^{-1} x(z\1 - a_2)^{-1} a_2 \de z  \\
&=& \frac{1}{2\pi i} \oint_\Gamma ((\lambda + z - z)\1 - a_1) ((\lambda + z\1) - a_1)^{-1} 
x (z\1 - a_2)^{-1} \de z  \\
&& + \frac{1}{2\pi i} \oint_\Gamma ((\lambda + z\1) - a_1)^{-1} 
x(z\1 - a_2)^{-1} (a_2 - z \1 + z \1)\de z  \\
&=& x\frac{1}{2\pi i} \oint_\Gamma (z\1 - a_2)^{-1} \de z  
+ \frac{1}{2\pi i} \oint_\Gamma ((\lambda + z\1) - a_1)^{-1} \de z \cdot x \\
&=& x\cdot \1 + \0 \cdot x = x. 
\end{eqnarray*}
This shows that $T$ is an inverse of $\lambda\1 - \Delta$, so that 
$\lambda\in{\mathbb C}\setminus\sigma_{\Lc(\Ac)}(\Delta)$. 
\end{proof}

\begin{corollary}\label{nil}
Let $\Ac$ be a complex FC-complete CIA. 
Assume that $a,b\in\Ac$ satisfy the condition 
$(\ad_{\Ac}a-\lambda)^mb=0$ for some $\lambda\in{\mathbb C}\setminus\{0\}$ 
and $m\ge1$, 
where $\ad_{\Ac}a\colon\Ac\to\Ac$, $(\ad_{\Ac}a)x=ax-xa$. 
Then for every integer  $N>2r_{\Ac}(a)/\vert\lambda\vert$ 
we have $b^N=0$. 
\end{corollary}

\begin{proof}
Since $\Delta:=\ad_{\Ac}a\colon\Ac\to\Ac$ is a derivation of $\Ac$, 
it follows by induction on $k$ that for all integers $n,k\ge1$, 
all $\lambda_1,\dots,\lambda_k\in{\mathbb C}$, 
and all $x_1,\dots,x_k\in\Ac$ 
we have 
$$(\Delta-\lambda_1-\cdots-\lambda_k)^n(x_1\cdots x_k)
=\sum_{j_1+\cdots+j_k=n}\frac{n!}{j_1!\cdots j_k!}
((\Delta-\lambda_1)^{j_1}x_1)\cdots((\Delta-\lambda_k)^{j_k}x_k).$$ 
Thence $(\Delta-k\lambda)^n(b^k)=0$ whenever $k\ge1$ 
and $n>km$. 

This shows that if $N\ge1$ and $b^N\ne0$, 
then $N\lambda\in\sigma_{\Lc(\Ac)}(\Delta)$. 
Consequently $N\lambda\in\sigma_{\Ac}(a)-\sigma_{\Ac}(a)$ by 
Proposition~\ref{rosenblum}, whence necessarily 
$N\vert\lambda\vert\le 2r_{\Ac}(a)$. 
\end{proof}

\section{Commutativity modulo the radical}

The following statement is a version of Lemma~1 in \S 24 of \cite{BS01}. 

\begin{lemma}\label{lemma1}
Let $\Ac$ be a unital complex FC-complete CIA with  $\rad\Ac=\{0\}$. 
Assume that $\g$ is a complex Lie subalgebra of $\Ac$ 
such that the closed unital associative subalgebra generated by $\g$ 
is equal to $\Ac$,  
and let $\fj$ be a finite-dimensional ideal of $\g$ 
such that $(\ad_{\g}a)|_{\fj}\colon\fj\to\fj$ 
is a nilpotent map for every $a\in\g$. 
Then $[\g,\fj]=\{0\}$. 
\end{lemma}

\begin{proof} 
We essentially 
follow the lines of the proof of Lemma~1 in \S 24 of \cite{BS01}, 
by using the previous Proposition~\ref{fine} instead 
of the Kleinecke-Sirokov Theorem. 
Since $\dim\fj<\infty$, it follows from Lemma~\ref{lemma5} that 
there exists some $m\ge1$ such that 
\begin{equation}\label{m}
(\forall a_1,\dots,a_m\in\g)(\forall a_0\in\fj)\quad 
(\ad_{\g}a_m)\cdots(\ad_{\g}{a_1})a_0=0.
\end{equation}
We shall prove that if $m\ge 2$ then \eqref{m} also holds 
with $m-1$ instead of $m$,  
and this will conclude the proof. 

To this end let $a_1,\dots,a_{m-1}\in\g$ and $a_0\in\fj$ arbitrary, 
and denote $y=(\ad_{\g}a_{m-1})\cdots(\ad_{\g}a_1)a_0\in\fj$. 
We have to check that $y=0$. 
By \eqref{m} we have 
$[a_m,y]=(\ad_{\g}a_m)(\ad_{\g}a_{m-1})\cdots(\ad_{\g}{a_1})a_0=0$ 
for all $a_m\in\g$. 
Since $\Ac$ is generated by $\g$ it then follows that 
$ay=ya$ for all $a\in\Ac$, that is, $y\in\Zc_{\Ac}$. 

On the other hand, note that $y=[a_{m-1},u]$, 
where $u=(\ad_{\g}a_{m-2})\cdots(\ad_{\g}a_1)a_0$ 
if $m\ge 3$ and $u=a_0$ if $m=2$. 
Since $y$ commutes with every element in $\Ac$, 
we have in particular $ya_{m-1}=a_{m-1}y$ and $yu=uy$, 
hence Proposition~\ref{fine} shows that $\sigma(y)=\{0\}$. 
Consequently, by Lemma~\ref{center}, we get 
$y\in\Qc_{\Ac}\cap\Zc_{\Ac}\subseteq\rad\Ac=\{0\}$, 
and the proof ends. 
\end{proof}

Here is a suitable version of Proposition~1 in \S 24 of \cite{BS01}. 

\begin{proposition}\label{proposition1}
Let $\Ac$ be a complex FC-complete CIA. 
Assume that $\g$ is a complex Lie subalgebra of $\Ac$ 
such that the closed unital associative subalgebra generated by $\g$ 
equals $\Ac$,  
and let $\fj$ be a finite-dimensional solvable ideal of $\g$. 
Then $\Nc_{\Ac}\cap\fj$ is an ideal of $\g$ and 
$\Nc_{\Ac}\cap\fj\subseteq\rad\Ac$. 
\end{proposition}

\begin{proof} 
We follow the lines of the proof of Proposition~1 in \S 24 of \cite{BS01}. 
Thus, let $\rho\colon\Ac\to\Lc(\Ac)$ be the regular representation of $\Ac$. 
Then $\rho(\g)$ is a Lie subalgebra of $\Lc(\Ac)$ 
and $\rho(\fj)$ is a finite-dimensional solvable ideal of $\rho(\g)$. 

On the other hand, let $\h$ be a Cartan subalgebra of $\fj$, 
$(\fj^\alpha)_{\alpha\in R}$ the root spaces of 
$\fj$ corresponding to the set $R$ of non-zero roots, 
and $\fj=\h\dotplus\g_+$ the corresponding Fitting decomposition, 
as in Proposition~\ref{theorem2} and its proof. 
It follows by Corollary~\ref{nil} that for every root $\alpha\in R$ 
we have $\fj^\alpha\subseteq\Nc_{\Ac}$, 
and hence $\rho(\fj^\alpha)\subseteq\Nc_{\Lc(\Ac)}$. 
Now note that $\rho\colon\Ac\to\Lc(\Ac)$ is a faithful representation, 
hence $\rho(\h)$ is a Cartan subalgebra of $\rho(\fj)$ and 
$(\rho(\fj^\alpha))_{\alpha\in R}$ are the corresponding root spaces. 
Also, $\rho(\Nc_{\Ac})\subseteq\Nc_{\Lc(\Ac)}$. 
Thus, we can apply Proposition~\ref{theorem2} to deduce that 
$\rho(\Nc_{\Ac}\cap\fj)$ is an ideal of $\rho(\fj)$. 
Since $\rho$ is faithful, this shows that 
$\Nc_{\Ac}\cap\fj$ is an ideal of $\fj$. 

Now, to prove that $\Nc_{\Ac}\cap\fj$ is even an ideal of $\g$, 
it suffices to check that $[a,\Nc_{\Ac}\cap\fj]\subseteq\Nc_{\Ac}$ 
for arbitrary $a\in\g$. 
To this end, denote $\fj_1={\mathbb C}a+\fj$, 
which is a finite-dimensional solvable Lie subalgebra of $\g$ 
since $\fj$ is a finite-dimensional solvable ideal. 
Then the preceding argument shows that $\Nc_{\Ac}\cap\fj_1$ 
is an ideal of $\fj_1$. 
If it happens that $\Nc_{\Ac}\cap\fj_1=\Nc_{\Ac}\cap\fj$, 
then $[a,\Nc_{\Ac}\cap\fj]=[a,\Nc_{\Ac}\cap\fj_1]\subseteq 
\Nc_{\Ac}\cap\fj_1=\Nc_{\Ac}\cap\fj$ and we are done. 
Now assume that $\Nc_{\Ac}\cap\fj_1\ne\Nc_{\Ac}\cap\fj$ 
and pick $b\in(\Nc_{\Ac}\cap\fj_1)\setminus\fj$. 
Then $\fj_1={\mathbb C}b+\fj=(\Nc_{\Ac}\cap\fj_1)+\fj$, hence 
$$\dim((\Nc_{\Ac}\cap\fj_1)/(\Nc_{\Ac}\cap\fj))
=\dim((\Nc_{\Ac}\cap\fj_1)/((\Nc_{\Ac}\cap\fj_1)\cap\fj))
=\dim(((\Nc_{\Ac}\cap\fj_1)+\fj)/\fj)
=\dim(\fj_1/\fj)\le1.$$
Since $b\in(\Nc_{\Ac}\cap\fj_1)\setminus(\Nc_{\Ac}\cap\fj)$, 
it then follows that 
$\Nc_{\Ac}\cap\fj_1={\mathbb C}b+(\Nc_{\Ac}\cap\fj)$. 
Now the finite-dimensional Lie algebra 
$\Nc_{\Ac}\cap\fj_1$ consists of nilpotent elements, hence it is nilpotent. 
And every hyperplane subalgebra of a finite-dimensional nilpotent Lie algebra 
is an ideal, 
so that $\Nc_{\Ac}\cap\fj$ is an ideal of $\Nc_{\Ac}\cap\fj_1$. 
Since $\Nc_{\Ac}\cap\fj$ is an ideal of $\fj$ it then follows that 
$$[a,\Nc_{\Ac}\cap\fj]\subseteq[\fj+(\Nc_{\Ac}\cap\fj_1),\Nc_{\Ac}\cap\fj]
\subseteq(\Nc_{\Ac}\cap\fj)+(\Nc_{\Ac}\cap\fj)=\Nc_{\Ac}\cap\fj.$$
Since $a\in\g$ is arbitrary, this completes the proof of the fact that 
$\Nc_{\Ac}\cap\fj$ is an ideal of $\g$.

It remains to show that $\Nc_{\Ac}\cap\fj\subseteq\rad\Ac$. 
By Lemma~\ref{lemma5}, it follows that there exists an integer $m\ge1$ 
such that the product of any $m$ elements in $\Nc_{\Ac}\cap\fj$ vanishes. 
On the other hand, we have already seen that $\Nc_{\Ac}\cap\fj$ is an ideal 
of $\g$, that is, 
$[\Nc_{\Ac}\cap\fj,\g]\subseteq\g$. 
This implies that 
$(\Nc_{\Ac}\cap\fj)\cdot\g
\subseteq\g\cdot(\Nc_{\Ac}\cap\fj)+(\Nc_{\Ac}\cap\fj)$, 
hence the unital associative subalgebra $\Ac_0$ of $\Ac$ 
generated by $\g$ satisfies 
$$(\Nc_{\Ac}\cap\fj)\cdot\Ac_0=\Ac_0\cdot(\Nc_{\Ac}\cap\fj),$$ 
where for any subsets $S_1,\dots,S_k$ of $\Ac$ 
we denote by $S_1\cdots S_k$ the linear subspace 
generated by all products $s_1\cdots s_k$ 
with $s_1\in S_1,\dots,s_k\in S_k$. 
By iterating the above equality $m$ times we get 
$$(\Nc_{\Ac}\cap\fj)\cdot\Ac_0\cdots(\Nc_{\Ac}\cap\fj)\cdot\Ac_0
=(\Nc_{\Ac}\cap\fj)\cdots(\Nc_{\Ac}\cap\fj)\cdot\Ac_0\cdots\Ac_0
=\{0\}.$$
In particular, for all $c\in\Nc_{\Ac}\cap\fj$ and $d\in\Ac_0$ 
we have $(cd)^m=0$. 
Since $\Ac_0$ is dense in $\Ac$ by hypothesis, 
the latter equality actually holds for all $d\in\Ac$ and implies that 
$\1-cd\in\Ac^\times$ for all $d\in\Ac$. 
That is, $c\in\rad\Ac$ for arbitrary $c\in\Nc_{\Ac}\cap\fj$, 
and the proof ends. 
\end{proof}

Now we can extend the main result of \cite{Ti87} 
(or Theorem~1 in \S 24 of \cite{BS01}) to 
FC-complete locally convex algebras with continuous inversion.

\begin{theorem}\label{turovskii}
Let $\Ac$ be a complex FC-complete CIA. 
Assume that $\g$ is a complex Lie subalgebra of $\Ac$ 
and let $\Ac(\g)$ be the closed unital associative subalgebra of $\Ac$ generated 
by $\g$. 
Then for every finite-dimensional solvable ideal $\fj$ of $\g$ 
we have $[\fj,\g]\subseteq\rad(\Ac(\g))\subseteq\Qc_{\Ac}$.  
\end{theorem}

\begin{proof}
The inclusion $\rad(\Ac(\g))\subseteq\Qc_{\Ac(\g)}$ 
($\subseteq\Qc_{\Ac}$) follows at once 
from the definition of $\rad(\Ac(\g))$. 

To prove $[\fj,\g]\subseteq\rad(\Ac(\g))$, 
first note that according to Lemma~\ref{subalgebras}(1), $\Ac(\g)$ is 
in turn a unital complex FC-complete CIA. 
Thus we may (and do) assume that $\Ac(\g)=\Ac$. 
Then denote $\widetilde{\Ac}=\Ac/\rad\Ac$ and let 
$q\colon\Ac\to\widetilde{\Ac}$ the natural projection. 
We claim that $\tilde\Ac^\times = q(\Ac^\times)$. The inclusion 
$q(\Ac^\times) \subeq \tilde \Ac^\times$ is trivial. 
For the converse, we assume that $q(a) \in \tilde\Ac^\times$. 
Then there exists $b \in \Ac$ such that $ab,ba \in \1 + \rad \Ac \subeq \Ac^\times$. 
Hence $a$ is left and right invertible, which implies that $a \in \Ac^\times$. 
Now $\1 + \rad \Ac \subeq \Ac^\times$ implies that 
$$q^{-1}(\tilde \Ac^\times) = q^{-1}(q(\Ac^\times)) = \Ac^\times + \rad \Ac 
= \Ac^\times,$$ 
and this entails that $\sigma_{\tilde\Ac}(q(a)) = \sigma_\Ac(a)$ for each $a \in A$. 
Now Remark~\ref{radical} and Lemma~\ref{subalgebras}(3) show   
that $\widetilde{\Ac}$ is 
a unital complex FC-complete CIA with 
$\rad\widetilde{\Ac}=\{0\}$. 
We are going to apply Lemma~\ref{lemma1} to 
the ideal $q(\fj)$ of the Lie subalgebra $q(\g)$ of $\widetilde{\Ac}$. 

In fact, let $a\in q(\g)$ arbitrary and $b\in q(\fj)$ such that 
$(\ad_{q(\g)}a)b=\lambda b$ for some $\lambda\in{\mathbb C}\setminus\{0\}$. 
Then 
$b\in\Nc_{\widetilde{\Ac}}\cap q(\fj)\subseteq\rad\widetilde{\Ac}=\{0\}$ 
by Corollary~\ref{nil} and Proposition~\ref{proposition1}. 
Consequently the linear mapping 
$(\ad_{q(\g)}a)|_{q(\fj)}\colon q(\fj)\to q(\fj)$ 
has no non-zero eigenvalue, and thus it has to be a nilpotent map. 
It then follows by Lemma~\ref{lemma1} 
that $[q(\g),q(\fj)]=\{0\}$, that is, 
$[\g,\fj]\subseteq\Ker q=\rad\Ac$. 
\end{proof}

\begin{corollary}\label{modulo}
Let $\Ac$ be a complex FC-complete CIA. 
Assume that $\g$ 
is a finite-dimensional complex Lie subalgebra of $\Ac$, 
$\fr$ its solvable radical,  
and let $\Ac(\g)$ be the closed unital associative subalgebra of $\Ac$ generated 
by $\g$. 
Then $[\g,\fr]\subseteq\rad(\Ac(\g))\subseteq\Qc_{\Ac}$. 
\end{corollary}

\section{CIA-linear Lie groups} 

\begin{definition} A finite-dimensional Lie group $G$ is said to be 
{\it linear} if it is isomorphic to a (closed) Lie subgroup of some 
$\GL_n(\R)$. Obviously, this condition is equivalent to the requirement that 
$G$ is isomorphic to a Lie subgroup of some finite-dimensional unital  algebra 
$\Ac$.   

We call $G$ {\it CIA-linear} if there exists an injective continuous homomorphism 
$\eta \: G \to \Ac^\times$ for some FC-complete CIA $\Ac$. 
\end{definition}

\begin{lemma} \label{rel-comp} 
Let $\Ac$ be a complex FC-complete CIA and 
$x \in \Ac$ quasi-nilpotent such that $e^{\R x}$ is relatively compact in $\Ac^\times$. 
Then $x = 0$. 
\end{lemma}

\begin{proof} Let $\Ac_0 \subeq \Ac$ be a maximal commutative subalgebra containing 
$x$. Then $\Ac_0$ is an FC complete commutative CIA with 
$\{0\} = \sigma_{\Ac}(x) = \sigma_{\Ac_0}(x)$ (Lemma~\ref{subalgebras}). 
We may therefore assume that $\Ac$ is commutative.  

Then the Gelfand transform $\Gc \: \Ac \to C(\hat \Ac), a \mapsto \hat a$ satisfies 
$\|\hat a\|_\infty = r_{\Ac}(a)$, so that 
$\rad{\Ac} = \Qc_\Ac = \ker \Gc$ is a closed $2$-sided ideal of $\Ac$. 
We conclude that the closure 
$K \subeq \Ac^\times$ of $e^{\R x}$ is contained in the closed 
affine subspace $U := \1 + \rad {\Ac} = \Gc^{-1}(\1)$, which is 
contained in $\Ac^\times$, hence $K$ is a subgroup of $\Ac^\times$. 

Next we observe that the Spectral Mapping Theorem implies that 
$$ \Exp \: (\Qc_\Ac,+) = (\rad(\Ac),+) \to U, \quad x \mapsto e^x $$
is a diffeomorphism whose inverse is given by the logarithm function 
$\Log \: U \to \rad(\Ac)$, which in turn is given by a convergent power series. 
We conclude that $\Log(K)$ is a compact subgroup of the locally convex space 
$(\rad(\Ac),+)$, hence trivial, and this implies that $x =0$.   
\end{proof}

The main result of our paper is the following theorem: 

\begin{theorem}\label{final} 
For a connected finite-dimensional Lie group $G$ 
the following are equivalent: 
\begin{itemize}
\item[(1)] $G$ is CIA-linear.  
\item[(2)] The continuous homomorphisms $\eta \: G \to \Ac^\times$ 
into the unit groups of FC-complete CIAs separate the points of $G$. 
\item[(3)] $G$ is linear. 
\end{itemize}
\end{theorem}

\begin{proof} (1) $\Rarrow$ (2) and (3) $\Rarrow$ (1) hold trivially. 
Therefore it remains to show that (2) implies (3). 

Assume that $G$ satisfies (2). 
Let $\g := \L(G)$ be the Lie algebra of $G$ and 
$\g = \fr \rtimes \s$ be a Levi decomposition. 

{\bf Step 1:} For $0 \not= x \in [\fr,\g]$ the subgroup 
$\exp_G(\R x)$ is not relatively compact. 

We argue by contradiction. Let $0 \not= x \in [\fr,\g]$. 
From (2) and the fact that $\g$ is finite-dimensional 
it follows that there exists a complex FC-complete CIA $\Ac$ 
and a morphism of Lie groups 
$\eta \: G \to \Ac$ for which $\L(\eta) \: \g \to \Ac$ is injective. 
Now Corollary~\ref{modulo} implies that 
$\L(\eta)([\g,\fr]) \subeq \Qc_\Ac$, so that 
Lemma~\ref{rel-comp} entails that $\eta(\exp_G(\R x)) = e^{\R \L(\eta)x}$ 
is not relatively compact in $\Ac^\times$, and this implies that 
$\exp_G(\R x)$ cannot be relatively compact in $G$. 

{\bf Step 2:} Let $R := \la \exp_G\fr \ra$ denote the radical of $G$. 
Then $R$ is a linear Lie group. 

Let $R'$ be the commutator subgroup of $R$. 
Since $\L(R') = [\fr,\fr]$, 
Step 1 implies that $\exp_G(\R x)$ is not relatively compact
in $G$ for $0 \not = x \in [\fr,\fr]$. This implies that $R'$ is closed 
(\cite{Ho65}, XVI.2.3/4) and does not contain circle groups. Hence it is a 
simply connected nilpotent Lie group, so that all compact subgroups of $R'$ are trivial. 
Now \cite{Ho65}, XVIII.3.2 implies that $R$ is a linear Lie group. 

{\bf Step 3:} The Levi subgroup $S := \la \exp_G \s \ra \leq G$ is linear. 

Let $q_S \: \tilde S \to S$ denote the universal covering and 
$\eta_S\: \tilde S \to \tilde S_\C$ be the universal complexification. 
Then $S$ is linear if and only if $S$ is a quotient of 
$\eta_S(\tilde S)$, i.e., $\ker \eta_S \supeq \ker q_S$ (\cite{Ho65}, XVII.3.3). 
Since $\Ac$ is a complex FC-complete CIA, the homomorphism of Lie algebras 
$\L(\eta) \: \s \to \Ac$ extends to a homomorphism $\L(\eta)_\C \: \s_\C \to \Ac$, 
which in turn integrates to a morphism of Lie groups $\eta_\C \: \tilde S_\C \to \Ac^\times$ 
with $\L(\eta_\C) = \L(\eta)_\C$, which implies that 
$\L(\eta_\C \circ \eta_S)\res_{\s} = \L(\eta)\res_\s$, and hence that 
$\eta_\C \circ \eta_S = \eta \circ q_S$, so that 
$\ker \eta_S \subeq q_S^{-1}(\ker \eta)$. By (2), the homomorphisms 
$\eta \: G \to \Ac^\times$ separate the points of $S$, we conclude that 
$\ker \eta_S \subeq \ker q_S$, showing that $S$ is linear. 

{\bf Step 4:} $G$ is linear because $R$ and $S$ are linear 
(\cite{Ho65}, XVIII.4.2). 
\end{proof}

\textbf{Acknowledgment.}
The first-named author was partially  supported 
from grant CERES 4-187/2004.

\end{document}